\documentclass[12pt]{amsart}     
\usepackage{amssymb}     
     
\begin{document}     
     
\newtheorem{theorem}{Theorem}[section]     
\newtheorem{lemma}[theorem]{Lemma}     
\newtheorem{corollary}[theorem]{Corollary}     
\newtheorem{proposition}[theorem]{Proposition}     
     
\theoremstyle{definition}     
\newtheorem{remark}{Remark}     
\def\theremark{\unskip}     
\newtheorem{definition}{Definition}     
\def\thedefinition{\unskip}     
\newtheorem{problem}{Problem}[section]     
     
\numberwithin{equation}{section}     
     
\def\Re{\operatorname{Re\,} }     
\def\dist{\operatorname{dist\,} }     
\def\p{\partial}     
\def\rp{^{-1}}     
     
\def\gap{\smallskip\noindent}     
\def\mgap{\medskip\noindent}     
\def\hils{{\mathcal H}}     
\def\R{{\Bbb R}}     
\def\Z{{\Bbb Z}}     
\def\Q{{\Bbb Q}}     
\def\C{{\Bbb C}}     
     
\def\reals{ {{\mathbb R}} }     
\def\naturals{ {{\mathbb N}} }     
\def\integers{ {{\mathbb Z}} }     
\def\complex{{\mathbb C}}     
\newcommand\beq{\begin{equation}}     
\newcommand\eeq{\end{equation}}     
\def\scriptd{ {\mathcal D} }     
\def\scriptdplus{ {\mathcal D}_+ }     
\def\scriptdminus{ {\mathcal D}_- }     
\def\scriptdpm{ {\mathcal D}_\pm }     
\def\scriptdmp{ {\mathcal D}_\mp }     
\def\scriptm{{\mathcal M}}     
\def\phila{\phi_\lambda}     
\def\scripte{ {\mathcal E} }  
\def\scriptt{{\mathcal T}}     
\def\scriptf{{\mathcal F}}     
\def\scriptg{{\mathcal G}}     
\def\scriptv{{\mathcal V}}   
\def\sump{\sideset{}{'}\sum}     
\def\marrow{\overset{\rightharpoonup}{m}}     
\def\jarrow{\overset{\rightharpoonup}{\jmath}}     
     
\author{Michael Christ}     
\address{     
        Michael Christ\\     
        Department of Mathematics\\     
        University of California \\     
        Berkeley, CA 94720-3840, USA}     
\email{mchrist@math.berkeley.edu}     
\thanks{The first author was supported in part by NSF grant DMS-9970660 
and completed this research while on appointment as a 
Miller Research Professor in the Miller Institute for Basic Research 
in Science.}     
\author{Alexander Kiselev}     
\address{ Alexander Kiselev\\     
Department of Mathematics\\     
University of Chicago\\     
Chicago, Ill. 60637}     
\email{kiselev@math.uchicago.edu}     
\thanks{The second author was supported in part by NSF grant DMS-9801530}     
     
\date\today{}     
     
\title[Absolutely continuous spectrum of Stark operators]     
{Absolutely continuous spectrum \\ of Stark operators}     
     
\begin{abstract}     
We prove several new results on the absolutely continuous spectra of   
perturbed one-dimensional  
Stark operators. First, we find new classes of perturbations,  
characterized mainly by smoothness  
conditions, which preserve purely absolutely continuous spectrum. Then  
 we establish stability of   
the absolutely continuous spectrum in more general situations, where  
imbedded singular spectrum   
may occur. We present two kinds of optimal conditions for the stability of  
 absolutely continuous   
spectrum: decay and smoothness. In the decay direction, we show that a 
sufficient (in the power scale) condition is   
 $|q(x)| \leq C(1+|x|)^{-\frac{1}{4}-\epsilon};$ in the smoothness  
direction, a sufficient condition in H\"older classes is   
$q \in C^{\frac{1}{2}+\epsilon}(\reals)$.   
On the other hand, we show that there exist potentials which both satisfy  
$|q(x)| \leq C(1+|x|)^{-\frac14}$   
and  belong to $C^{\frac12}(\reals)$ for  
which the spectrum becomes purely singular  
on the whole real axis, so that the above results are optimal within the  
scales considered.  
\end{abstract}     
     
\maketitle     
     
\section{Introduction}     
    
In this paper we consider the Stark operator   
\begin{equation}  
\label{stark}   
H_q = -\frac{d^2}{dx^2} -x +q(x)   
\end{equation}  
defined on the whole real line $\reals.$ This operator describes a charged quantum  
particle in a constant electric field  subject to an  
additional electric potential $q(x).$ There exists an extensive physical and   
mathematical literature on Stark operators; for a review, see e.g. \cite{CFKS}.  
When $q(x)=0,$ the operator has purely absolutely continuous   
spectrum. The question we wish to address is which classes of perturbations  
 $q$ preserve this property. We will consider two classes of conditions  
that ensure preservation of the absolutely continuous spectrum:   
smoothness and decay.  
The first result on the smoothness condition was proven by Walter   
\cite{Wal}, who showed that if the potential is bounded and   
has two bounded derivatives,   
the spectrum remains purely absolutely continuous. Similar results were  
obtained by Bentosela, Carmona, Duclos, Simon, Souillard and Weder  
in \cite{Ben} using Mourre's method. A corollary noted in \cite{Ben}  
is a drastic change in the spectral properties of Schr\"odinger   
operators of Anderson model type, say   
\[ H_\omega = -\frac{d^2}{dx^2} + \sum\limits_n a_n(\omega) V(x-n) \]  
where $a_n$ are independent identically distributed random variables  
and $V \in C_0^2 ((0,1)),$ when a constant electric field is switched on.  
The spectrum changes from almost surely pure point to purely absolutely   
continuous.   
Recently, Sahbani  
\cite{Sa1,Sa2} relaxed smoothness  
conditions of \cite{Wal} and \cite{Ben} (see the remark after Theorem~\ref{spac2}).  
On the opposite side of the smoothness scale, Delyon, Simon and Souillard \cite{DSS}   
showed that for a periodic array of $\delta$ function potentials with random couplings  
in a constant electric field, the spectrum is purely singular.   
Avron, Exner and Last  
\cite{AEL} realized that the spectrum may  be purely singular even   
for a deterministic periodic array of very singular interactions, such  
 as $\delta'.$  
Generalizations of these results, as well as other models with singular   
potentials, were considered in \cite{Min,Ex,ADE,MS,Briet,ABD}.   
There remained a  gap, however, between the classes of   
potentials for which localization was known to occur, and those
 for which the   
spectrum was known to remain absolutely continuous.   
  
As far as decay conditions are concerned, it is well known that if   
$q(x)$ satisfies $|q(x)| \leq C(1+|x|)^{-\alpha},$ $\alpha>1/2,$ then the    
spectrum remains purely absolutely continuous \cite{Va}.  
Moreover, there are examples where $|q(x)x^{1/2}|\leq C$ and   
isolated imbedded eigenvalues appear. If   
$|q(x)|x^{1/2} \rightarrow \infty,$ it was shown by Naboko and Pushnitski  
\cite{NaPu} that  
 dense (imbedded) point spectrum   
may appear on all of $\reals$.  We remark that for the   
operator without electric field, the decay threshold where imbedded  
eigenvalues may appear is the power $-1;$ of course, it is physically natural  
that it is more difficult to get an imbedded eigenvalue in the presence  
of the constant electric field.  
However, if we do not wish to rule out imbedded singular spectrum, 
it has  been shown in \cite{Ki1} that the absolutely continuous spectrum   
of a perturbed Stark operator still  fills the whole real axis 
when $|q(x)| \leq C(1+|x|)^{-\alpha},$ $\alpha>1/3.$   
The question what is the critical rate of decay for which   
the spectrum may become purely singular remained open.   
  
Our main goal in this paper is to prove two sharp results on the preservation  
of the absolutely continuous spectrum of Stark operators.   
 Recall that $f(x)$ is   
called H\"older continuous with exponent   
$\alpha$ ($f \in C^\alpha(\reals)$)  
if  
\[ \|f\|_{C^\alpha} = {\rm sup}_x |f(x)| + {\rm sup}_{x,y}  
\frac{|f(x)-f(y)|}{|x-y|^\alpha} < \infty. \]  
 
We will analyze solutions of the generalized eigenfunction equation 
\begin{equation}  \label{starke} 
-u'' -xu + q(x)u = Eu\ . 
\end{equation} 
Here $-x$ represents a background potential due to a constant 
electrical field, while $q$ is some perturbation. 
  
\begin{theorem}  
\label{sac1}  
Assume that the potential $q(x)$ is H\"older continuous with exponent $\alpha>1/2.$  
Then an essential support of the absolutely continuous part of the spectral  
measure coincides with the whole real axis. Moreover, for a.e.~$E,$ all  
solutions $u(x,E)$ of equation   
\eqref{starke} satisfy $u(x,E)=O(x^{-1/4}),$ $u'(x,E)=O(x^{1/4})$ as $x \rightarrow +\infty.$  
\end{theorem}  
  
\it Remark. \rm 1. An essential support of $\mu$ is a set $S$  
such that $\mu(\reals \setminus S)=0$ and  
$\mu (S_1)>0$ for any $S_1 \subset S$ of positive Lebesgue measure. \\  
2. In this and subsequent theorems, only behavior of $q(x)$  
for $|x|$ large matters. We will always implicitly assume  
$q$ to be locally integrable, and will state only additional  
hypotheses which concern its behavior for large $x$. 
On the negative part of the real axis, it is sufficient  
for all our conclusions 
to require that $q(x)-x \rightarrow +\infty$ as $x \rightarrow -\infty$.   
We prefer to state the results in a slightly weaker form to avoid   
making statements too cumbersome.   
  
\begin{theorem}  
\label{dac1}  
Assume that the potential $q(x)$ is locally integrable, 
and that $q(x^2) \in L^p$ for some $1 \leq p <2.$  
Then an essential support of the absolutely continuous part of the  
spectral measure coincides with the whole real axis.   
Moreover, for a.e.~$E$,  all solutions $u(x,E)$ of equation   
\eqref{starke} satisfy $u(x,E)=O(x^{-1/4}),$ $u'(x,E)=O(x^{1/4})$  
as $x \rightarrow +\infty.$    
\end{theorem}  
\it Remark. \rm 1. In particular,  
the assumption of Theorem~\ref{dac1} is satisfied if 
$|q(x)| \leq C(1+|x|)^{-\alpha}$, for some $\alpha>1/4$.  
\\ 
2. An explicit expression for the leading term in an asymptotic 
expansion of $u(x,E)$, as $x\to+\infty$, 
can also be derived under the hypotheses of 
Theorems~\ref{sac1} and \ref{dac1};   
see Theorem~\ref{genthm} and Section~\ref{sm}.  
  
Both Theorems~\ref{sac1} and \ref{dac1} are direct corollaries of the  
following more general result.   
\begin{theorem}\label{genthm}  
Consider a Stark operator $H_q$ on $\reals^1$. 
Assume that the potential $q(x)$ admits a decomposition   
$q=q_1+q_2,$ where both $q_1(x^2)$ and $x^{-1}q'_2(x^2)$ belong   
to $(L^1+L^p)(\reals),$ and that there exists $\zeta<1$ such that   
$|q_2(x)|\leq \zeta |x|$ for sufficiently large $|x|$.  
Then for almost every     
energy $E$ there exists a solution $u_+(x,E)$ of equation  
\eqref{starke} with the asymptotic behavior  
\begin{equation}  
\label{wkbas}  
u_+(x,E) = (x-q_2(x)+E)^{-1/4}   
e^{i\phi(x,E)} 
\left(1+o(1)\right)  
\end{equation}  
 as $x \rightarrow +\infty$, 
where    
\begin{equation*} 
\phi(x,E) 
= \int\limits_0^x 
\Big[\sqrt{t-q_2(t)+E}-\frac{q_1(t)}{2\sqrt{x-q_2(t)+E}}\Big] 
\,dt 
\end{equation*} 
\end{theorem}  
\it Remark. \rm   
A sufficient condition on the derivative is that  
$|q'(x)| \leq C(1+|x|)^{\alpha}$, for some $\alpha < 1/4$.  
A surprising aspect of this theorem is that the perturbation $q_2$ 
is allowed  virtually as much growth as the constant electric field  
potential, and  more flexibility on the derivative.  
  
Theorem~\ref{dac1} follows immediately; Theorem~\ref{sac1} requires a  
 simple argument showing that   
any $C^\alpha$ potential with $\alpha>1/2$ can be represented as in  
 Theorem~\ref{genthm}.   
We sketch this argument in Section~\ref{sm}.   
  
We notice that in the case of a Schr\"odinger operator without constant   
electric field, the absolutely continuous spectrum is preserved for potentials  
with power decay rate $\alpha>1/2$ \cite{CK,Re,DK}.   
There exist potentials $V$ satisfying $|V(x)x^{1/2}| \leq C$ for which   
the absolutely continuous spectrum is destroyed \cite{KU,KLS},  
so that $\alpha = 1/2$ is a sharp threshold.  
It is natural that in the presence of a constant electrical field, 
the absolutely continuous spectrum is preserved under more slowly decaying 
perturbations of the potential. 
  
The next result shows optimality of Theorems~\ref{sac1} and \ref{dac1}.  
Fix $f \in C_0^\infty ((0,1)),$ not identically zero,   
and let $a_n(\omega)$ be independent,   
identically distributed random variables with uniform distribution in   
$[0,2\pi].$ Denote    
$c=\left({3}/{2}\right)^{2/3}.$ Let us define  
\begin{equation}  
\label{qx}  
q(x) = c \sum\limits_{n=1}^\infty n^{-1/2} f(\sqrt{cx}-n)  
\sin (\tfrac{4}{3} x^{\frac{3}{2}} + a_n(\omega)) \ . 
\end{equation}  
(There is nothing magic in the choice of $c;$ however, this choice  
will simplify computations later.)  
 We have   
\begin{theorem}  
\label{count}  
Let $q(x)$ be a random potential given by \eqref{qx}. Then for a.e.~$\omega,$   
the spectrum of the corresponding perturbed Stark operator   
is purely singular on the whole real line.  
\end{theorem}  
  
In particular, any realization of $q$ defined by \eqref{qx}  
satisfies $|q(x)| \leq Cx^{-1/4}$ and   
 belongs to $C^{1/2}(\reals),$ so Theorem~\ref{count}  
assures sharpness of Theorems~\ref{sac1} and \ref{dac1} in H\"older spaces  
and in the power scale, respectively.   
  
On a more detailed level, we may want to distinguish between  
 perturbations for which the absolutely continuous  
 spectrum is preserved, but imbedded   
singular spectrum may appear, and perturbations  
  which preserve purely absolutely continuous spectrum.   
Our final results provide two criteria ensuring pure absolute   
continuity of the spectrum.  
 Recall that a function $g(x)$ defined   
on a real line is called smooth in Zygmund's sense if   
\[ \int\limits_0^1 {\rm sup}_{x \in \reals} |g(x+\epsilon)-2g(x)+  
g(x-\epsilon)| \frac{d \epsilon}{\epsilon^2} < \infty. \]  
\begin{theorem}  
\label{spac1}  
Assume that the potential $q(x)$ is bounded,   
has bounded continuous first derivative and is smooth in Zygmund's sense.  
Then the Stark operator \eqref{stark} has purely absolutely continuous   
spectrum on the whole real axis.   
\end{theorem}  
  
We can also allow the potential to grow at a rate arbitrarily close   
to that of the constant electric field, and still have purely absolutely   
continuous spectrum, provided that we impose a slightly different condition   
on smoothness.  
  
\begin{theorem}  
\label{spac2}  
Assume that the potential $q(x)$ satisfies $q(x) = O(x^\alpha)$  
as $|x|\to\infty$ for some   
$\alpha<1,$ is differentiable, and that its derivative $q'(x)$ is Dini   
continuous:  
\[ \int\limits_0^1 {\rm sup}_x |q'(x+\epsilon)-q'(x-\epsilon)|  
\frac{d \epsilon}{\epsilon} < \infty. \]  
Then the spectrum of the perturbed Stark operator \eqref{stark} is purely   
absolutely continuous on the whole real axis.   
\end{theorem}  
\it Remark. \rm The classes of potentials that 
are smooth in Zygmund's  sense or have Dini continuous derivative  
were first considered in this context by Sahbani in \cite{Sa1, Sa2}.  
Using the conjugate operator approach,   
he proved that under conditions similar to Theorem~\ref{spac1} or   
Theorem~\ref{spac2} (but with stronger growth restrictions) the spectrum   
is absolutely continuous with perhaps some imbedded eigenvalues.   
His results extend (in a slightly weaker form) to the higher dimensional  
setting. \\  
  
We will later discuss examples demonstrating that 
the criteria given by Theorems~\ref{spac1} and \ref{spac2}  
are fairly sharp.   
  
We employ two different approaches to prove the stated results.   
First,  
to prove  Theorems~\ref{spac1}, \ref{spac2}, and \ref{count},  
we apply a Liouville transformation to reduce the Stark operator   
to a form reminiscent of the Schr\"odinger operator without electric  
field, but with the energy entering in a non-standard way.  
We then use a Pr\"ufer transformation  
to analyze the asymptotic behavior of solutions.  
In the proofs of the other results, it is more convenient to  
represent \eqref{starke} as a first-order system and to employ 
estimates for the solution series. Some of the tools for these estimates  
come from our recent work \cite{CK,CK1,CK2,CKN}.   
As soon as we have control over the asymptotic behavior   
of solutions, we can apply the whole axis version of subordinacy theory   
due to Gilbert \cite{Gi}, or the approximate eigenvectors criterion of \cite{CKL},  
to draw spectral conclusions.   
  
\section{Preservation of purely absolutely continuous spectrum}  
  
We begin by proving Theorems~\ref{spac1} and \ref{spac2} as a warm-up.   
All the proofs of spectral properties in this paper rely on the   
study of solutions of the equation \eqref{starke}, 
$-u'' -xu +q(x)u = Eu$. 
The link between   
the behavior of solutions and spectral results is provided by   
Gilbert-Pearson subordinacy   
theory, more particularly, by the whole-line version   
of this theory due to Gilbert \cite{Gi}. Recall that a real solution   
$u_1(x,E)$ of \eqref{starke} is called subordinate on the right  
if for any other real linearly independent solution $u_2(x,E)$ we have   
\[ \lim_{N \rightarrow \infty} \frac{\int\limits_0^N |u_1(x,E)|^2\,dx}  
{\int\limits_0^N |u_2(x,E)|^2 \,dx} =0. \]  
Subordinacy on the left is defined similarly. 
Note that it is easy to see that for equation \eqref{starke}, under the 
assumptions of any of our theorems, there is always a solution subordinate  
(in fact, $L^2$) on the left since the potential goes to $+\infty$ there.   
The main result of Gilbert implies that singular spectrum may only be   
supported on the set of energies where there exists a solution  
subordinate  on  both sides. Moreover, the set of the energies  
where there exists a solution   
subordinate on one side, but there is no subordinate solution on the other  
side, is an essential support of the absolutely continuous spectrum, of   
multiplicity one. Therefore, our goal is to prove that for all  
energies (if we want to show pure absolute continuity), or for a.e. energy  
(if we allow imbedded singular spectrum), there is no solution of   
\eqref{starke} subordinate on the right.   
  
In the equation \eqref{starke}, let us perform a Liouville transformation  
given by (see, e.g. \cite{Ol})  
\begin{equation}  
\label{liou}  
 \xi(x) = \int\limits_0^x \sqrt{t} \, dt = \tfrac{2}{3}x^{3/2},  
\,\,\,   
\phi (\xi) = x(\xi)^{1/4} u(x(\xi)).   
\end{equation}  
This transformation introduces an irrelevant singularity at the origin; 
henceforth we always work outside some neighborhood of $0$. 
 
The resulting function $\phi$ satisfies the Schr\"odinger equation   
\begin{equation}  
\label{basic1}  
-\phi'' + \left( \frac{5}{36 \xi^2} + \frac{q(c\xi^{2/3}) -E}  
{c \xi^{2/3}} \right) \phi = \phi,   
\end{equation}  
where $c=({3}/{2})^{2/3}$. Let us introduce   
a short-hand notation $V(\xi,E)$ for the expression in brackets in   
\eqref{basic1}. Let us further apply a Pr\"ufer transformation   
to the equation for $\phi,$ setting for each $E$  
\begin{equation}  \label{proofer} 
\begin{split}  
 \phi(\xi,E) & = R(\xi,E) \sin (\theta(\xi,E)) \\  
 \phi'(\xi,E) & = R(\xi,E) \cos (\theta(\xi,E)).  
\end{split}  
\end{equation} 
The equations for  $R$ and $\theta$ are as follows:   
\begin{align}  
\label{pru1}  
 (\log R(\xi, E))'  
& =   
\tfrac{1}{2}V(\xi,E)\sin (2\theta(\xi,E)) \\  
\label{pru2}  
\theta'(\xi,E)  
& =   
1- \tfrac{1}{2}V(\xi,E)(1-\cos (2\theta(\xi,E))).   
\end{align}  
Our main goal in proof of Theorem~\ref{spac1} will   
be to show the convergence of the integral   
\begin{equation}  
\label{keyint}  
 \int\limits_{1}^{N} \left( \frac{5}{36\xi^{2}} +   
\frac{-E + q(c\xi^{\frac{2}{3}})}{c\xi^{\frac{2}{3}}}\right)  
\sin (2\theta(\xi,E))\, d \xi   
\end{equation}  
as $N \rightarrow \infty$ for every $E.$ This goal is motivated by  
\begin{proposition}  
\label{GP}  
Suppose that   
\[ \limsup_{x \rightarrow \infty}|q(x)|/x = \zeta <1, \]  
 and for a given $E$ the integral   
\eqref{keyint} converges for all initial values of $\theta(0,E).$   
Then for this value of $E,$ there is no subordinate on the right solution   
of the equation \eqref{starke}.   
\end{proposition}  
\begin{proof}  
 If for a given value of $E$ the integral \eqref{keyint}  
 converges, it follows from \eqref{pru1} and \eqref{pru2} that all solutions   
of the equation \eqref{basic1} are bounded  and, moreover,   
any solution $\phi_\beta$ (where $\beta$ parametrizes  
boundary condition) has the following asymptotic behavior   
as $\xi \rightarrow +\infty:$  
\begin{equation*}  
\phi_{\beta}(\xi,E) = C_\beta \sin (\xi + g_{\beta}(\xi,E))(1+o(1)),   
\end{equation*}  
where $|g_{\beta}'(\xi,E)|<\zeta<\zeta_1<1$ for $\xi$ sufficiently large.  
Going back to the original equation \eqref{starke},   
we infer that   
every solution $u_{\beta}(x,E)$   
 has the following asymptotic behavior as   
$x \rightarrow \infty:$  
\[ u_{\beta}( x, E) =  
C_\beta x^{-\frac{1}{4}}  
\sin (\tfrac{2}{3}x^{\frac{3}{2}}+ f_{\beta}(x,E))(1+o(1)), \]  
where $|f_{\beta}'(x,E)|<\zeta_1 x^{1/2}$ for sufficiently large $x.$   
For any $u_{\beta}(x,E)$ we find  
\begin{multline} \label{normg}  
\int\limits_{1}^{N} |u_{\beta}(x,E)|^{2}\,dx = C_\beta^2\int\limits_{1}^{N}  
x^{-\frac{1}{2}}(\sin(\frac{2}{3}x^{\frac{3}{2}}+ f_{\beta}(x,E))^{2}\,   
dx \left( 1+o(1)\right)  \\  
= C_\beta^2  
 \left( N^{\frac{1}{2}} - \frac{1}{2}\int\limits_{1}^{N}x^{-\frac{1}{2}}  
\cos (\tfrac{4}{3}x^{\frac{3}{2}}+2 f_{\beta}(x,E))\, dx \right)  
\left(1+o(1)\right).   
\end{multline}  
Consider the integral  
\begin{equation}  
\label{auxint}  
I_\beta(N)= \int\limits_{1}^{N}x^{-\frac{1}{2}}  
e^{\left(i\frac{4}{3}x^{\frac{3}{2}}+2 f_{\beta}(x,E)\right)}\, dx.   
\end{equation}  
Integrating by parts, with $e^{2i f_\beta(x,E)}$ being differentiated,   
we obtain that \eqref{auxint} is equal to   
\begin{equation*}  
 -\int\limits_N^\infty x^{-\frac{1}{2}}  
e^{i\frac{4}{3}x^{\frac{3}{2}}} \, dx e^{2i f_\beta(N,E)}   
+ 
2i \int\limits_1^N f'_\beta(x,E)  
e^{2i f_\beta(x,E)} \int\limits_x^\infty t^{-\frac{1}{2}}  
e^{i\frac{4}{3}t^{\frac{3}{2}}} \, dt dx.  
\end{equation*}  
Since   
\[  \int\limits_x^\infty t^{-\frac{1}{2}}  
e^{i\frac{4}{3}t^{\frac{3}{2}}} \, dt = \frac{1}{2}x^{-1}e^{i\frac{4}{3}  
x^{\frac{3}{2}}}(1+o(1)), \]  
we obtain that   
\[ |I_\beta (N)| \leq 2\zeta_1 N^{1/2} \]  
for $N$ sufficiently large. Returning to \eqref{normg}, it is straightforward  
to conclude that any solution $u_\beta(x,E)$ satisfies for sufficiently   
large $N$   
\[ C_{\beta}^2(1-\zeta_1)  
 N^{1/2} \leq \int\limits_{1}^{N} |u_{\beta}(x,E)|^{2}\,dx   
\leq C_\beta^2 (1+\zeta_1) N^{1/2}. \]    
Therefore, all solutions have the same rate of $L^2$ norm growth  
as $N \rightarrow \infty,$ and there is no subordinate solution.   
\end{proof}   
  
Now we establish convergence of \eqref{keyint} for every energy   
under the assumptions of   
Theorem~\ref{spac1}, thus completing the proof of this result.  
\begin{proof}[Proof of Theorem~\ref{spac1}]  
Let us write  
\[ \theta_\beta (\xi,E) = \xi +g_\beta(\xi,E), \]  
with $\beta$ parametrizing the initial condition at $0$ and   
\begin{equation}  
\label{control}  
|g'_\beta(\xi,E)| \leq C\xi^{-2/3}  
\end{equation}  
 uniformly in $\beta.$  
Clearly we can ignore the short-range quadratic decay term and consider  
only   
\[ \int\limits_1^N \frac{E-q(c\xi^{2/3})}{c\xi^{2/3}} e^{2i(\xi +  
g_\alpha(\xi,E))}\,d\xi. \]  
Moreover, the integral  
\begin{multline*}  
 \int\limits_1^N \xi^{-2/3} e^{2i(\xi +  
g_\alpha(\xi,E))}\,d\xi  
\\ 
= -\left( \int\limits_\xi^\infty   
\frac{e^{2i\eta}}{\eta^{2/3}}\,d\eta \right) e^{2ig_\alpha(\xi,E)}\left|_0^N  
\right. 
+  \int\limits_1^N e^{2ig_\alpha(\xi,E)} g'_\alpha(\xi,E)  
\int\limits_\xi^\infty \frac{e^{2i\eta}}{\eta^{2/3}}\,d\eta d\xi  
\end{multline*}  
is clearly convergent due to \eqref{control}. It remains to estimate  
\[ \int\limits_{3}^{N}  \frac{q(c\xi^{2/3})}{c\xi^{2/3}} e^{2i(\xi +  
g_\alpha(\xi,E))}\,d\xi \]  
uniformly as $N \rightarrow \infty$   
(we shifted the region of integration for convenience).  
Fix $h,$ $\pi/2 > h > 0,$ and consider the equality  
\begin{multline*}  
\left(e^{2ih}+e^{-2ih}-2 \right)   
\int\limits_{3}^{N}  \frac{q(c\xi^{2/3})}{c\xi^{2/3}} e^{2i(\xi +  
g_\alpha(\xi,E))}\,d\xi  
= O(1) +  
\\   
\int\limits_{3}^{N}   
e^{2i\xi} \left[ \frac{q(c(\xi+h)^{2/3}) e^{2ig_\alpha  
(\xi+h,E)}}{(\xi+h)^{2/3}}+  
 \frac{q(c(\xi-h)^{2/3}) e^{2ig_\alpha(\xi-h,E)}}{(\xi-h)^{2/3}}-  
 \frac{2q(c\xi^{2/3}) e^{2ig_\alpha(\xi,E)}}{\xi^{2/3}} \right]  
\\  
= O(1)  
+ \int\limits_{3}^{N} \frac{e^{2i(\xi+g_\alpha(\xi,E))}}{\xi^{2/3}}  
\left[ q(c(\xi+h)^{2/3})+q(c(\xi-h)^{2/3}) -2 q(c \xi^{2/3}) \right]\, d\xi.  
\end{multline*}  
Since we assumed that $q \in C^1,$ it suffices to control  
\begin{equation}  
\label{con2}  
 \int\limits_{3}^{N} \frac{1}{\xi^{2/3}}   
\left| q(c(\xi^{2/3} + \tfrac{2}{3} \xi^{-1/3} h)) +   
q(c(\xi^{2/3} + \frac{2}{3} \xi^{-1/3} h))- 2q(c\xi^{2/3})\right|\,d\xi.   
\end{equation}  
Set $\epsilon = \frac{2}{3} c \xi^{-1/3}h,$ then  
 uniformly in $N,$ the integral \eqref{con2} is bounded by   
\[ C \int\limits_0^1 \left|q(c_1 \epsilon^{-2} +\epsilon) +   
q(c_1 \epsilon^{-2} -\epsilon)- 2 q(c_1 \epsilon^{-2})\right|  
\frac{d \epsilon}{\epsilon^2} \]  
which is finite by assumption.   
\end{proof}  
\it Remark. \rm Without change, the proof goes through   
even with the weaker growth assumption $q(x) = O(x^{\frac{1}{2}-\epsilon})$   
for some $\epsilon>0.$  
  
\begin{proof}[Proof of Theorem~\ref{spac2}]  
The proof of this theorem is very similar to the preceding proof.   
However, it is convenient to employ a slight variation of the Pr\"ufer  
transformation. Namely, we let  
\begin{eqnarray*}  
\sqrt{1-V(\xi,E)}\phi (\xi,E) = \tilde{R} \sin(\tilde{\theta}(\xi,E)) \\  
\phi'(\xi,E) = \tilde{R} \cos(\tilde{\theta}(\xi,E)).   
\end{eqnarray*}  
This transformation is well-defined for large $\xi$ where $V(\xi,E)<1,$   
and this suffices for our purpose since we are interested in the asymptotic  
behavior at $+\infty.$ The equations for $\tilde{R}$ and $\tilde{\theta}$   
are   
\begin{align}  
\label{tilR}  
(\log \tilde{R})'(\xi,E)  
& = -\frac{V'(\xi, E)}{4(1-V(\xi,E))}(1-  
\cos 2\tilde{\theta} (\xi,E)), \\  
\label{tilthe}  
\tilde{\theta}'(\xi,E)  
& = \sqrt{1-V(\xi,E)} -  
\frac{V'(\xi, E)}{4(1-V(\xi,E))}\cos 2\tilde{\theta} (\xi,E).   
\end{align}   
The role analogous to the integral \eqref{keyint} is played by   
\begin{equation}  
\label{keyint1}  
\int\limits_1^N \frac{V'(\xi,E)}{1-V(\xi,E)} \cos 2\tilde{\theta}(x,E) \,d\xi;  
\end{equation}  
the other term on the right hand side of \eqref{tilR} can be integrated explicitly.  
If the integral \eqref{keyint1} converges for a given energy $E$, then for this energy   
there is no solution subordinate on the right. This can be shown   
in a direct analogy to the proof of Proposition~\ref{GP}; the details are  
left to the reader.   
  
Expressing $V'/(1-V)$ in terms of $q,$   
we see that it is enough to show the convergence of  
$ \int_{a}^{N} q'(c \xi^{2/3}) \xi^{-1} e^{i \tilde{\theta}(\xi,E)}\,d\xi $ 
as $N \rightarrow \infty.$   
From \eqref{tilthe} and the assumption of the theorem, it follows that   
\[ \tilde{\theta}(\xi,E) = \xi + \tilde{g}_\alpha(\xi,E), \]  
where $|\tilde{g}'_\alpha(\xi,E)| \leq C \xi^{-\delta}$ for some  
$\delta>0.$   
Now fix $h,$ $0<h<\pi/2,$ and consider  
\begin{multline*}  
(e^{ih}-e^{-ih}) \int\limits_a^N    
q'(c \xi^{2/3}) \xi^{-1} e^{i \tilde{\theta}(\xi,E)}\,d\xi  
\\ 
=  O(1) +  
\int\limits_a^N e^{i\xi} \left[ e^{i\tilde{g}_\alpha (\xi+h)}  
\frac{ q'(c(\xi+h)^{2/3})}{\xi+h}-  e^{i\tilde{g}_\alpha (\xi-h)}  
\frac{ q'(c(\xi-h)^{2/3})}{\xi-h} \right]\, d\xi  
\\ 
=   
O(1) + \int\limits_a^N e^{i\xi+i\tilde{g}_\alpha(\xi)} \xi^{-1}   
\left[  q'(c(\xi+h)^{2/3})- q'(c(\xi-h)^{2/3})\,d\xi \right]  
\\ 
\leq 
C\left (1+  \int\limits_a^N \xi^{-1}   
|q'(c(\xi+h)^{2/3})- q'(c(\xi-h)^{2/3})| \,d\xi \right).  
\end{multline*}  
Setting   
\[ \epsilon = \frac{c(\xi+h)^{2/3} -c(\xi-h)^{2/3}}{2}, \]  
and making a change of variable in the last integral, we find that for   
a sufficiently large $a,$ the controlling integral \eqref{keyint1} is bounded   
by   
\[ C\left(1+ \int\limits_0^1 |q'(f(\epsilon)+\epsilon)-q'(f(\epsilon)  
-\epsilon)| \frac{d \epsilon}{\epsilon}\right) \]  
which is finite  
by assumption of Dini continuity  
(here $f(\epsilon) =  [c(\xi+h)^{2/3}+c(\xi-h)^{2/3}]/{2}$).  
\end{proof}  
  
We remark that the results of Theorems~\ref{spac1} and \ref{spac2} are  
rather sharp. For example, let $E=0,$ then \eqref{starke} reduces  
to   
\[ -\phi'' +\left( \frac{5}{36\xi^2}+  
\frac{q(c\xi^{2/3})}{c\xi^{2/3}}\right) \phi = \phi. \]  
Let us denote by $wn(\xi)$ the classical Wigner-von Neumann potential  
\cite{WvN, ReSi}. Choose $q$ so that the expression in the brackets   
coincides with $wn(\xi)$ for $\xi \geq 1.$ It is not difficult to show   
that we can take $q=0$ on $(-\infty,0)$ and $q$ smooth and bounded on $(0,1)$   
so that the whole equation \eqref{starke} has an eigenvalue at $E=0.$  
(The issue is gluing together the $L^2$ solution on $-\infty$ and the   
$L^2$ solution produced by Neumann-von Wigner potential on $\infty.$  
It can always be achieved by choosing $q$ appropriately on $(0,1):$  
see, e.g. \cite{Sim2} for a similar argument.) Notice that in this case  
\[ q(x)= x\left( wn(\tfrac{3}{2}x^{3/2})-\frac{5}{81x^3}\right). \]  
The Wigner-von Neumann potential has asymptotic behavior \cite{ReSi}  
\[ wn(x) = -8(\sin 2x)/x+O(x^{-2}), \]  
with the $O(x^{-2})$ term also smooth with derivatives decaying at   
the same rate. We see that   
\[ q(x) = C_1 x^{-1/2} \sin (C_2 x^{3/2})+q_1(x), \]  
where $q_1(x)$ is better behaved in all respects. This function $q(x)$   
narrowly misses the class of functions smooth in Zygmund's sense.   
Indeed, $|q(x+\epsilon)+q(x-\epsilon)-2q(x)| = \epsilon q_2(x,\epsilon),$   
where $q_2$ is uniformly bounded.   
Also, $q$ has a bounded,  
continuous derivative, which however fails to be Dini continuous.    
  
\section{Main Theorem}  
  
Here we prove Theorem~\ref{genthm}.   
As before, we are going   
to study the asymptotic behavior of the solutions to equation \eqref{starke},  
\[ -u''-xu(x)+q(x)u(x)=Eu(x), \]  
as $x \rightarrow \infty.$  The $L^1$ part of the perturbation can be treated  
by standard means (such as, for example, Levinson's theorem), so we will  
assume that $q_1(x^2) \in L^p(0,\infty)$ and $x^{-1} q'_2(x^2) \in   
L^p(0,\infty).$    
Write \eqref{starke} as a system  
\begin{equation}  
\label{sys1}  
 \left( \begin{array}{c} u \\ u' \end{array} \right)' = \left(  
\begin{array}{cc} 0 & 1 \\ -x+q(x)-E & 0 \end{array} \right) \left(  
 \begin{array}{c} u \\   
u' \end{array} \right).   
\end{equation}  
We are going to perform a series of transformations with this system,  
similarly to \cite{CK1,CKN}.  
Applying first a variation of parameters-type transformation   
\begin{equation}  
\label{tran1}  
\left( \begin{array}{c} u \\ u' \end{array} \right) =   
\left(  
\begin{array}{cc} e^{i\psi(x,E)} &   
e^{-i\psi(x,E)}  
 \\ i\psi'(x,E) e^{i \psi(x,E)} &   
-i \psi'(x,E) e^{-i \psi(x,E)}  
\end{array} \right) z  
\end{equation}  
we arrive at   
\[ z' = \left(  
\begin{array}{cc} -i\scripte & -i\overline{\scripte}e^{-2i\psi} \\  
 i\scripte e^{2i\psi}  & i\overline{\scripte}  \end{array} \right) z, \]  
where   
\begin{equation}\label{errdef}  
 \scripte(x,E) = \frac{1}{2\psi'(x,E)}\left(-i\psi''(x,E) + \psi'(x,E)^2  
 -x +q(x)-E \right).   
\end{equation}  
Letting   
\begin{equation}\label{tran2}  
 z= \left( \begin{array}{cc} e^{-i\int_0^x \scripte (t,E)\,dt } & 0 \\ 0 &  
   e^{i\int_0^x \overline{\scripte} (t,E)\,dt }  
\end{array} \right)y   
\end{equation}  
leads to   
\begin{equation}  
\label{findec}  
y' = \left( \begin{array}{cc} 0 & -i\overline{\scripte}e^{-2i\psi +i\int_0^x  
\Re \scripte (t,E)\,dt} \\  
 i\scripte e^{2i\psi -i\int_0^x \Re \scripte(t,E)\,dt} &   
0 \end{array} \right) y.   
\end{equation}  
We are going to choose   
\[ \psi'(x,E)= \sqrt{x-q_2(x)+E}, \]  
so that   
\begin{equation}\label{err}  
\scripte (x,E) = \frac{q_1(x)}{2\sqrt{x-q_2(x)+E}}- i\frac{1-q'_2(x)}{2(x-q_2(x)+E)}.  
\end{equation}  
Let $Q(x) = (x^{-1}q_1(x)^2+x^{-2}q_2'(x)^2+x^{-2})^{1/2},$   
and   
\begin{equation}\label{adef}  
 a(x,E) =  i \scripte(x,E) Q(x)^{-1}.   
\end{equation}  
Let us introduce multilinear operators  
\begin{multline}     
\label{muop}     
S_n (f_1, \dots,f_n)(x,E)   \\   =   \int\limits_x^\infty     
\int\limits_{t_{1}}^\infty     
\cdots    
\int\limits_{t_{n-1}}^\infty     
\prod_{j=1}^n    
\Big[    
e^{2i(-1)^{n-j}\left(\psi(t_j,E) -\int\limits_0^{t_j} \Re   
\scripte (t,E) \right)\,dt}a_j(t_j,E)  
f_j(t_j) dt_j    
\Big]    
\end{multline}  
and   
\begin{multline}       
\tilde{S}_n (f_1, \dots,f_n)(\xi,E)   \\   = \int\limits_x^\infty     
\int\limits_{t_{1}}^\infty     
\cdots    
\int\limits_{t_{n-1}}^\infty     
\prod_{j=1}^n    
\Big[    
e^{2i(-1)^{n-j}\left(\psi(t^2_j,E) -\int\limits_0^{t^2_j} \Re  
\scripte (t,E) \right)\,dt}a_j(t_j^2,E)  
f_j(t_j) dt_j    
  \Big],    
\end{multline}  
where $a_j$ is equal to $a$ for $n-j$ even and $\overline{a}$ for $n-j$ odd.    
Set $\tilde{f}(t) = 2f(t^2)t.$ Notice that making in \eqref{muop} a  
 change of variables $t_j=s^2_j$ leads to  
\begin{equation}  
\label{rel1}  
S_n (f_1, \dots, f_n)(x , E) = \tilde{S} (\tilde{f}_1,\dots,   
\tilde{f}_n)(x^{1/2},E).  
\end{equation}   
Notice that by assumption,   
\[ \tilde{Q}(x)=2(q_1(x^2)^2+x^{-2}q_2'(x^2)^2+x^{-2})^{1/2} \]  
belongs to $L^p$ with $p<2.$  
Iterating  system \eqref{findec} starting from the vector      
$(1,0)$ and using \eqref{rel1}, we obtain the following     
 formal series expansion for one of     
 the solutions:     
\begin{equation}     
\label{keyser}       
y_+(x,E) =  
 \left( \begin{array}{c}    
\sum\limits_{n=0}^{\infty}     
 S_{2n}(Q, \dots ,Q  
 )(x,E)  \\      
-\sum\limits_{n=1}^{\infty}     
S_{2n-1}(Q, \dots ,Q )(x,E)     
\end{array} \right)=  
 \left( \begin{array}{c}    
\sum\limits_{n=0}^{\infty}     
 \tilde{S}_{2n}(\tilde{Q}, \dots , \tilde{Q})(x^{1/2},E)  \\      
-\sum\limits_{n=1}^{\infty}     
 \tilde{S}_{2n-1}(\tilde{Q}, \dots, \tilde{Q})(x^{1/2},E)     
\end{array} \right)  
\end{equation}    
(we stipulate $S_0(Q)(\xi,E)=  
\tilde{S}_0(\tilde{Q})(\xi,E)\equiv 1$ in the above formula).      
Introduce an operator   
\begin{equation}  
\label{basop}  
\tilde{S}f(E) = \int\limits_0^\infty   
 \tilde{S}(t,E)  
f(t) dt,    
\end{equation}  
where   
\[ \tilde{S}(t,E) = a(t^2,E)  
e^{2i \left(\psi(t^2,E)   
- \int_0^{t^2}\Re \scripte (t,E)\right) \,dt}. \]  
  
\begin{proposition}  
\label{series}  
Fix a compact interval $J \subset \reals.$   
Assume that the operator   
$\tilde{S}$ maps $L^p(\reals)$ boundedly to  $L^r(J)$,  
for some $p <2<r.$ Then for any   
$\tilde{Q} \in L^p(\reals)$, the series \eqref{keyser} converges for   
a.e.~$E \in J.$ Moreover, for a.e.~$E \in J,$    
the solution $y_+(x,E)$ of the   
system \eqref{findec} given by \eqref{keyser} has the asymptotic behavior  
\begin{equation}  
\label{keyasym}  
y_+(x,E) =  \left( \begin{array}{c} 1 \\ 0 \end{array} \right) +o(1).   
\end{equation}  
\end{proposition}  
\begin{proof}  
The proof is based on results of \cite{CK1,CK2,CKN}.  
Introduce a multilinear operator $M_n$, acting on $n$ functions  
$g_k (x,E)$, by  
\[  
M_n(g_1,\dots,g_n)(x,x',E) = \int_{x\le x_1\le\cdots\le x_n\le x'}  
\prod_{k=1}^n \big[ g_k(x_k,E)\,dx_k\big]  
\ .  
\]  
In the special case when there is a single function $g$ such that  
each $g_k$ is either $g$ or $\overline{g}$, we write  
simply  
$M_n(g)(x,x',E)$.  
In particular, the multilinear transforms  $\tilde{S}$ in the series  
\eqref{keyser} have a structure  
identical to $M_n,$ with $g(x,E) = \tilde{S}(x,E)\tilde{Q}(x).$   
Define  
\begin{align*}  
M_n^*(g_1,\dots\,g_n)(E)   
& = \sup_{x\le x'\in\reals} \left| M_n(g_1,\dots\,g_n)(x,x',E) \right|  
\\  
M_n^*(g)(E)   
& =  
\sup_{x\le x'\in\reals} \left| M_n(g)(x,x',E) \right|  
\ ,  
\\  
G(g,E)  
& = \sum_{r=0}^1 \sum_{m=1}^\infty  
m \left( \sum_{\jmath=1}^{2^m} \Big|   
\int_{E^m_\jmath} g(x,E)\,dx  
\Big|^2 \right)^{1/2}  
\ .  
\end{align*}  
In Proposition~4.2 and in the proof of Theorem~1.3  
 of \cite{CK2} it is shown that  
\begin{align}  
M_n^*(g_1,\dots,g_n)(E)  
&\le C^n\prod_{k=1}^n G(g_k,E)  
\label{M_nbound1}  
\\  
M_n^*(g)(E)  
& \le C^n\frac{ G(g,E)^n}{\sqrt{n!}}  
\label{M_nbound2}  
\end{align}  
for some universal constant $C<\infty$.  
If the operator $\tilde{S}$ satisfies the required   
$L^p(\reals) - L^r(J)$ bound, then it is not hard to see that   
\begin{equation}  
\label{Gestimate}  
 \|\tilde{G}(\tilde{S}(Q)),E)\|_{L^r(J)} \leq   
C \|Q\|_{L^p(\reals)}   
\end{equation}  
(see Proposition~3.3 of \cite{CK1}).  
The estimates \eqref{M_nbound1}, \eqref{M_nbound2},   
\eqref{Gestimate} allow to show a.e.$~E$ convergence of the  
series \eqref{keyser}, and to prove Proposition~\ref{series}.  
For details we refer to \cite{CK1}, Section 4, where a similar argument   
is given.  
\end{proof}  
\it Remark. \rm An alternative route to the same result is to consider   
an energy-dependent potential $\scripte (x,E),$ rather than to introduce $a(x,E).$   
The paper \cite{CKN} follows this approach. \\  
  
It remains to show that the operator $\tilde{S}$ satisfies   
$L^p(\reals)-L^r(J)$ bound for any compact $J$ and   
any $1 \leq p <2,$ $r = p/(p-1)$  
provided that   
potential $q$ satisfies assumptions of Theorem~\ref{dac1}.  
Such result would follow by complex interpolation if we establish   
$L^1(\reals) \mapsto L^\infty(J)$ and $L^2(\reals) \mapsto L^2(J)$ bounds.  
The first bound is evident since $a(x^2,E)$ is bounded and the   
oscillatory exponential  
\[ e^{2i \left(\psi(x,E)- \int_0^{x}  \Re \scripte(t,E) \right)\,dt} \]  
is bounded too.   
Next we establish the key $L^2-L^2$ bound for the operator $\tilde{S}.$   
  
\begin{proposition}  
\label{l2b}  
Let $J \subset \reals$ be a compact interval.   
Assume that   
\begin{equation}  
\label{assa}  
|\partial_E^\beta a(x,E)| \leq C  
\end{equation}  
for $\beta=0,1,2$ for every $E \in J,$ and    
that the  potential  
$q$ satisfies the assumptions of Theorem~\ref{genthm}.   
Then for any   
$f \in L^2(\reals)$ we have   
\[ \|\tilde{S} f \|_{L^2(J)} \leq C \|f\|_{L^2(\reals)} \]  
for an operator $\tilde{S}$ defined by \eqref{basop}.  
\end{proposition}   
\begin{proof}  
Notice that   
\[ \|\tilde{S} f \|^2_{L^2(J)} \leq  
\int\limits_\reals \int\limits_\reals f(x) \overline{f}(y)   
\left(\int\limits_J  \eta(E)a(x^2,E) \overline{a}(y^2,E)  
e^{2i \int\limits_{x^2}^{y^2} \left(   
\psi'(t,E)- \Re \scripte(t,E) \right)\,dt}  
  \, dE \right) \,dxdy, \]  
where $\eta(E)$ is a positive $C_0^\infty(\reals)$ function satisfying   
$\eta(E) \geq 1$ for $E \in J.$   
We can rewrite the kernel in the brackets in the above formula by   
making a change of variable $t=s^2:$  
\begin{equation}  
\label{kernel}  
K(x,y) = \int\limits_J \eta(E)a(x^2,E) \overline{a}(y^2,E)  
e^{ 4i  \int\limits_y^x \left(   
s\psi'(E,s^2)-s \Re \scripte (s^2,E) \right)\,ds}  
  \, dE.  
\end{equation}  
Let us denote for simplicity   
\[ \sigma(s,E) = (\psi'(s^2,E)-\Re \scripte (s^2,E))s. \]  
Since   
\[ \partial_E (s\sqrt{s^2-q_2(s^2)+E})= \frac{s}{2\sqrt{s^2-q_2(s^2)+E}}  
=1+o(1) \]  
for large $s,$ a direct computation using the assumption on $q$ shows that   
 there exists $N$ such that for   
 any $x,y>0$ such that  $|x-y|\geq N,$ we have   
uniformly in $E \in J$   
\begin{equation}\label{phest}  
\left| \int\limits_x^y \partial_E \sigma(s,E) \,ds \right|   
\geq C_1 |x-y|, \,\,\,   
 \left| \int\limits_x^y \partial_E^\beta \sigma(s,E) \,ds \right|  
\leq C_2|x-y|, \,\,\,\beta=2,3    
\end{equation}  
with positive constants $C_1,$ $C_2.$   
Integrate by parts two times   
in \eqref{kernel}, with    
\[ e^{4i \int_{x}^{y} \sigma (s,E) \,ds }   
\int_{x}^{y} \partial_E \sigma (s,E) \, ds\] being integrated.  
Using \eqref{phest} and  
\eqref{assa}  
we obtain   
\[ |K(x,y)| \leq C (1 + |x-y|^2)^{-1} \]  
for all $x$ and $y,$ which implies the desired $L^2-L^2$ bound.  
\end{proof}  
The fact that $a(x,E)$ satisfies \eqref{assa} can be checked directly from definitions   
\eqref{errdef} and \eqref{adef}.  
  
Now we complete the proof of Theorem~\ref{dac1}.  
  
\begin{proof}[Proof of Theorem~\ref{dac1}]  
Proposition~\ref{l2b} ensures that   
 Proposition~\ref{series} applies under conditions of the theorem.  
Notice that   
\[ e^{i\psi(x,E)-i \int_0^x \scripte (t,E)\,dt}=  
 (x-q_2+E)^{-1/4} e^{i \int\limits_0^{x}  
\left( \sqrt{t-q_2(t)+E}-\frac{q_1(t)}{2\sqrt{t-q_2(t)+E}} \right)  
\,dt }. \]  
Applying transformations \eqref{tran2}, \eqref{tran1} and   
\eqref{liou} to the solution   
$y_+(\xi,E)$ and using \eqref{keyasym}, we obtain for a.e.~$E \in J$ a  
solution $u_+(x,E)$ of the equation \eqref{starke} with the   
asymptotic behavior   
\[ u_+(x,E) = (x-q_2+E)^{-1/4} e^{i \int\limits_0^{x}  
\left( \sqrt{t-q_2(t)+E}-\frac{q_1(t)}{2\sqrt{t-q_2(t)+E}} \right)  
\,dt }\left(1+o(1) \right), \]  
coinciding with \eqref{wkbas}, and  
\[ u'_+(x,E) = i(x-q_2+E)^{1/4} e^{i \int\limits_0^{x}  
\left( \sqrt{t-q_2(t)+E}-\frac{q_1(t)}{2\sqrt{t-q_2(t)+E}} \right)  
\,dt }\left(1+o(1) \right). \]  
There also exists solution $u_-(x,E)$ which is just complex conjugate   
of $u_+.$ In particular, for a.e.~$E,$ any solution $u(x,E)$ of \eqref{starke}  
satisfies   
\begin{equation}\label{decsolest}  
|u(x,E)| \leq Cx^{-1/4}, \,\,\, |u'(x,E)| \leq Cx^{1/4}.   
\end{equation}   
By the results of \cite{CKL}, to show that the essential support of the absolutely continuous part   
of the spectral measure fills all real  
line, it suffices to show that for a.e.~$E$ there exists a sequence $\psi_n(x,E)$ such that   
\begin{equation}\label{appcri}  
\limsup_{n \rightarrow \infty} \frac{|\psi_n(0,E)|+  
|\psi'_n(0,E)|}{\|\psi_n\| \| (H_q-E) \psi_n\|}>0.  
\end{equation}  
We can obtain an appropriate sequence $\psi_n(E)$ by taking a smooth cutoff function   
$\eta(x),$ $1>\eta(x)>0,$ $\eta(x) =1$ if $x<1/2,$ $\eta(x)=0$ if $x>1,$ and letting   
$\psi_n(x,E) = \tilde{u}(x,E) \eta(2^{-n}x),$ where $\tilde{u}$ is the solution of \eqref{starke}  
which is $L^2$ at $-\infty.$ Then $|\psi_n(0,E)|+|\psi'(0,E)|>c(E)>0$ since $\tilde{u}$   
is nonzero. Moreover, using \eqref{decsolest} gives  
$\|\psi_n\| \sim 2^{n/4}$ and $\|(H_q-E)\psi_n\| \sim 2^{-n/4},$ so that  
\[ \|\psi_n\| \| (H_q-E) \psi_n\| \leq C(E).\]  
Therefore, \eqref{appcri} is satisfied at a.e.~$E.$  
  
An alternative (but more technical) way to finish the proof is to apply subordinacy theory   
using more detailed information on the asymptotic behavior, similarly   
to Proposition~\ref{GP}.  
\end{proof}  
  
\section{Potentials in $C^{\frac{1}{2}+\epsilon}$}\label{sm}   
  
In this section we are going to prove a result from which Theorem~\ref{sac1}   
will follow. Define   
\[ \overline{D}^{\alpha}f(x) = {\rm sup}_{y:\,|x-y|<1}   
\frac{|f(x)-f(y)|}{|x-y|^\alpha}, \]  
$0<\alpha \leq 1.$  
\begin{theorem}\label{sacs}  
Assume that potential $q(x)$ satisfies   
$|q(x)| \leq \zeta |x|$ with $\zeta <1$ for all sufficiently   
large $|x|,$ and   
 $x^{-\alpha} \overline{D}^{\alpha}q(x^2) \in L^p(1,\infty)$ for   
some $p<2$ and $1 \geq \alpha>0.$   
Then  an essential support of the absolutely continuous part  
of the spectral measure coincides with the whole real axis.   
\end{theorem}  
\noindent \it Remark. \rm In particular,   
if $q \in C^\alpha(\reals)$ with $\alpha>1/2,$ then   
$D^\alpha q \in L^\infty$   
and the assumption of Theorem~\ref{sacs} is satisfied.  
\begin{proof}  
We will show that if $q$   
satisfies the conditions of the theorem, then it can be represented  
as a sum $q=q_1+q_2$ with $q_1(x^2) \in L^p(\reals)$ and $x^{-1}q'_2(x^2) \in L^p(\reals),$   
$|q_2(x)| \leq \zeta_1 x$ with $\zeta_1  <1$ for large $x.$   
Then the result follows from Theorem~\ref{genthm}.   
Fix a function $\eta(x) \in C_0^\infty(0,1)$ such that $\int_\reals \eta(x) \,dx =1$.   
Set   
\begin{align*}   
q_2(x) &= x^{1/2} \int_\reals \eta(x^{1/2}(x-y))q(y)\,dy, 
\\  
q_1(x) &= x^{1/2} \int_\reals \eta(x^{1/2}(x-y))(q(x)-q(y))\,dy.  
\end{align*}  
Then   
\[ |q_1(x)| \leq    
x^{(1-\alpha)/2} \int_\reals \eta(x^{1/2}(x-y))\overline{D}^{\alpha}q(x)\,dy   
= x^{-\alpha/2}\overline{D}^{\alpha}q(x). \]  
Therefore, by assumption, $q_1(x^2) \in L^p.$   
The property $|q_2(x)| \leq \zeta_1 x$ with $\zeta_1<1$   
for all large $x$ is clear from   
the definition.  
Also   
\begin{eqnarray}   
q'_2(x)&  = & \frac{1}{2} x^{-1/2}   
\int_\reals \eta(x^{1/2}(x-y))q(y)\,dy  \nonumber \\  
& & + x \int_\reals \eta'(x^{1/2}(x-y))   
(q(y)-q(x))\,dy + \frac{1}{2} \int_\reals   
\eta'(x^{1/2}(x-y))(x-y)q(y)\,dy. \label{q_2est}  
\end{eqnarray}  
Due to the bound $|q(x)| \leq \zeta x,$ the first term on the right hand side of \eqref{q_2est}  
is bounded, and so harmless.   
In the third term,  $|x-y| \leq x^{-1/2}$  where the integrand is nonzero, hence it is also   
bounded (using $x^{1/2}\int |\eta'(x^{1/2}(x-y)| \,dy \leq C$).  
Finally, in the second term we estimate   
\[ |q(x)-q(y)| \leq C x^{-\alpha/2} \overline{D}^\alpha q(x) \]  
where $\eta$ is nonzero, and so this term is bounded by   
$Cx^{(1-\alpha)/2}\overline{D}^\alpha q(x).$   
Hence, in all, $x^{-1}q'_2(x^2) \in L^p.$   
\end{proof}  
  
\section{A Counterexample}  
  
Our main goal in this section is to prove Theorem~\ref{count}, 
showing that  
the results of Theorems~\ref{sac1} and \ref{dac1} are sharp.   
Namely, we will prove that there exist potentials which both belong to   
$C^{1/2}(\reals)$, and decay at the rate $x^{-1/4}$, yet which   
lead to singular spectrum on the whole real axis. Hence one example  
will show optimality of both theorems.   
Apart from establishing optimality,  
the example we are going to construct  
is interesting in its own right, suggesting the mechanism   
to get singular spectrum in a background constant electric field  
``at the lowest cost'',  
and giving potentials with interesting dynamical properties.     
Our potentials will be random,  
and the construction is inspired   
by \cite{KLS}, although there are some notable differences.  
In \cite{KLS}, random perturbations of the free operator of the form   
\begin{equation}  
\label{KLSR}  
 V(x) = \sum\limits_{n=1}^\infty a_n(\omega)n^{-\alpha} f(x-n)   
\end{equation}  
were considered (here $f$ and $a_n$ can be taken as in \eqref{qx}).   
This model has a transition of spectral properties from absolutely   
continuous to singular at $\alpha = 1/2$ (see also \cite{KU} for a similar 
earlier model).   
  
By subordinacy theory, it is enough to construct the potential on the   
positive half-axis with the   
above decay and smoothness properties such that for a.e.~$E \in \reals$  
there exists a solution subordinate on the right. We are going to analyze  
the equation \eqref{basic1}  
\[-\phi'' + \left( \frac{5}{36 \xi^2} + \frac{q(c\xi^{2/3}) -E}  
{c \xi^{2/3}} \right) \phi = \phi, \]  
and its Pr\"ufer variables representation \eqref{pru1}, \eqref{pru2} 
\begin{align*}  
 (\log R(\xi, E))'  
& =   
\frac{1}{2}  \left( \frac{5}{36 \xi^2} + \frac{q(c\xi^{2/3}) -E}  
{c \xi^{2/3}} \right)  \sin (2\theta(\xi,E)) \\  
\theta'(\xi,E) & =   
1- \frac{1}{2} \left( \frac{5}{36 \xi^2} + \frac{q(c\xi^{2/3}) -E}  
{c \xi^{2/3}} \right)(1-\cos (2\theta(\xi,E))).   
\end{align*}  
Let us denote by $\beta$ the boundary condition for $\theta$ at the   
origin.  
Let $a_n(\omega)$   
be independent random variables with uniform distribution in $[0,2\pi].$   
Fix a function $f \in C_0^\infty((0,1)),$ $f$ not identically zero.  
Consider a family of random   
potentials $q(x)$ chosen so that   
\begin{equation}  
\label{ranpot}   
\frac{q(c\xi^{2/3})}{c\xi^{2/3}} = \xi^{-2/3} \sum\limits_{n=1}^{\infty}  
n^{-1/2} f(\xi^{1/3}-n) \sin (2\xi +a_n(\omega)).  
\end{equation}  
Notice that in the coordinate $x$, according to the Liouville transformation,   
our potential $q(x)$ looks exactly like \eqref{qx}:  
\[ q(x) = c \sum\limits_{n=1}^\infty n^{-1/2} f(\sqrt{cx}-n)  
\sin (\tfrac{4}{3} x^{\frac{3}{2}} + a_n(\omega)). \]  
In particular, any such $q$ satisfies $|q(x)| \leq Cx^{-1/4},$ and   
any such $q$ belongs to $C^{1/2}(\reals).$    
We begin the proof of Theorem~\ref{count} with a series of auxiliary   
statements.  
The key, as before, will be the analysis of the asymptotic behavior   
of solutions. First, we are going to show that in order to apply subordinacy  
theory, it is sufficient to study the asymptotic behavior of $R(\xi,E).$  
\begin{lemma}  
\label{derivest}  
For any $q(x)$ satisfying $|q(x)| \leq Cx,$   
any energy $E$ and any boundary condition $\beta,$   
we have   
\[ \int\limits_1^L \phi^2(\xi,E,\beta)\,d\xi \leq   
 \int\limits_1^L R^2(\xi,E,\beta)\,d\xi \leq  
C_q(E) \int\limits_1^L \phi^2(\xi,E,\beta)\,d\xi. \]  
\end{lemma}  
\begin{proof}  
Notice that condition on $q$ ensures that $q(c\xi^{2/3})\xi^{-2/3}$   
is bounded. By definition of $R,$ it is enough to show that the   
$L^2$ norm of the derivative of $\phi$ on $[1,L]$  
is controlled by some constant  
$C_q(E)$ times the $L^2$ norm of $\phi.$  
The proof of the latter fact is a simple elliptic regularity type argument,   
and is well known. We refer, for example, to \cite{LS1}, proof of   
Theorem 2.1C, for a proof of completely parallel result.      
\end{proof}  
  
\begin{theorem}  
\label{rangrowth}  
Consider the Pr\"ufer variables equations 
\eqref{pru1}, \eqref{pru2},  with potential $q$ defined  
by \eqref{ranpot}.  
Then for any $E \in \reals$ and any boundary condition $\beta$ we have   
for a.e.~$\omega$   
\begin{equation}  
\label{lyapbeh}  
\lim\limits_{\xi \rightarrow \infty} \frac{\log R(\xi,E,\beta)}{\log \xi} =  
\Lambda (E),   
\end{equation}  
where   
\begin{equation}  
\label{lyapexp}  
\Lambda(E) = \frac{3\pi}{8} |\hat{f} (\frac{3E}{C})|^2  
\end{equation}  
(here $\hat{f}$ denotes the Fourier transform $\int \exp(i\xi x)f(x)\,dx$ 
of $f$).   
\end{theorem}  
\it Remark. \rm Notice that ``Lyapunov exponent'' $\Lambda(E)$ is   
positive everywhere except for perhaps a discrete set of points   
since the Fourier transform of $f$ is analytic. \\  
  
We begin with a sequence of auxiliary lemmas. 
Consider the unique solution of 
\eqref{pru2} satisfying some fixed boundary condition $\beta$.
First, notice that the integral   
\[ \int\limits_1^N  
\left( \frac{5}{36 \xi^2} -\frac{E}{c \xi^{2/3} } \right) \sin  
(2 \theta (\xi,E)) \, d \xi \]  
stays bounded as $N \rightarrow \infty$ for any $\omega$ by the   
argument in the proof of Theorem~\ref{spac1},   
and so in the equation \eqref{pru1} for $R(\xi,E)$ it is enough to 
consider  the integral  
\begin{equation}  
\label{conint}  
 \frac{1}{2}\int\limits_1^N \frac{q(c \xi^{2/3})}{c \xi^{2/3}} 
\sin (2\theta  (\xi, E))\, d\xi.   
\end{equation}  
First we are going to analyze $R(\xi,E)$ along a sequence $\xi=n^3,$   
as suggested by the form of the potential \eqref{ranpot}, since the parts  
with independent phases are supported on $[n^3, (n+1)^3].$    
Notice that   
\begin{multline}  
\label{devest}  
\Big| \int\limits_{n^3}^{\xi}  
 \frac{q(c \xi^{2/3})}{c \xi^{2/3}} \sin (2\theta  
(\xi, E))\, d\xi \Big|  
\\ 
=  \Big| n^{-1/2}  
\int\limits_{n^3}^{\xi} \xi^{-2/3} f(\xi^{1/3}-n)  
\sin (2\xi+ a_n(\omega)) \sin (2\theta(\xi,E)) \, d\xi \Big|  
\leq Cn^{-1/2}   
\end{multline}  
for $\xi \in [n^3, (n+1)^3].$  
Therefore, if we control $R(n^3,E),$ we will have sufficient control  
for all $\xi,$ in particular   
\begin{equation}  
\label{intsame}  
\lim\limits_{x \rightarrow \infty} \frac{R(x,E)}{\log x}  
 =  \frac{1}{3} \lim\limits_{n \rightarrow \infty} \frac{R(n^3,E)}{\log n}.   
\end{equation}  
Denote by $\theta_n$ the value of $\theta$ at the beginning of the   
$n$th interval, $\theta_n = \theta(n^3,E) = \theta(n^3,E,\beta,\omega)$.   
  
Define 
\begin{equation} \label{Fdefinition} 
F(\xi,E) = \tfrac{1}{2}\int\limits_{n^3}^{\xi}\frac{E}{c\eta^{2/3}}\, d\eta  
= \frac{3}{2c}E(\xi^{1/3}-n).     
\end{equation} 
\begin{lemma}   
\label{thetarep}  
We can represent $\theta (\xi,E)$ on $[n^3, (n+1)^3]$ in the following form:  
\begin{multline}  
\theta(\xi,E)  =  \theta_n + \xi + F(\xi,E) + \Theta(\xi,E)  
\\   
\label{theta}  
+ \tfrac12 {n^{-1/2}} \int\limits_{n^3}^\xi  \eta^{-2/3} f(\eta^{1/3}-n)  
\sin (2\eta+ a_n(\omega)) \cos(2\eta + 2F(\eta, E)+2\theta_n)\, d\eta,   
\end{multline}  
where  
\begin{equation*} 
|\Theta(\xi,E)| \leq C(E) n^{-1}.  
\end{equation*} 
\end{lemma}  
 
\begin{proof}  
Let us analyze the terms arising in the equation \eqref{pru2} for the angle   
$\theta.$ The proof reduces to showing that   
\begin{multline*}   
\Theta(\xi,E)  =  \frac{1}{2}\int\limits_{n^3}^\xi  
 \left[   
\left(   
\frac{5}{36 \eta^2}  
+\frac{E}{c\eta^{2/3}}\right) \cos(2\theta(\eta,E)) -  
\left(  \frac{5}{36 \eta^2}+\frac{q(c\eta^{2/3})}{c\eta^{2/3}} \right)   
\right.  
\\   
\left.  
+\   
\frac{q(c\eta^{2/3})}{c\eta^{2/3}}\left( \cos(2\theta(\eta,E))  
- \cos(2\eta + 2F(\eta, E)+2\theta_n)\right)   
\right]   
\, d\eta   
\end{multline*}   
satisfies the required estimate.   
The first term satisfies   
\begin{equation}  
\label{aaa1}  
 \Big| \int\limits_{n^3}^\xi \Big( \frac{5}{36 \eta^2}  
+\frac{E}{c\eta^{2/3} }\Big) \cos(2\theta(\eta,E)) \, d\eta \Big|   
\leq C n^{-2}   
\end{equation}  
by an integration by parts argument identical to that in the proof of   
Theorem~\ref{spac1}. The second term satisfies   
\[ \big| \int\limits_{n^3}^\xi \Big(   
 \frac{5}{36 \eta^2}+\frac{q(c\eta^{2/3})}{c\eta^{2/3}} \Big)\, d \eta   
\Big| \leq Cn^{-3}, \]  
also by a simple integration by parts estimate for the second summand.  
For the last term, we have   
\begin{multline} \label{longint}   
\int\limits_{n^3}^\xi   
\frac{q(c\eta^{2/3})}{c\eta^{2/3}}\left( \cos(2\theta(\eta,E))  
- \cos(2\eta + 2F(\eta, E)+2\theta_n)\right)\, d\eta  
\\ 
=    
2\int\limits_{n^3}^\xi \frac{f(\eta^{1/3}-n)}{n^{1/2} \eta^{2/3}}  
\sin(2\eta +a_n(\omega)) \sin(\theta(\eta,E)+\eta+F(\eta,E)+\theta_n)  
\\  
\times  
\sin(\theta(\eta,E) - \eta -F(\eta,E)-\theta_n)\,d\eta.   
\end{multline}  
Now, we have by \eqref{pru2}  
\begin{multline*}  
|\theta (\eta, E) - \eta - F(\eta,E)-\theta_n|  
\\ =   
\left| \frac{1}{2}  
\int\limits_{n^3}^\eta   
\left[ \left(  
\frac{q(ct^{2/3})}{ct^{2/3}}+  
\frac{5}{36t^2}\right)(1-\cos 2\theta(t,E)) -  
\right. \right.  
\left. \left.   
\frac{E}{ct^{2/3}} \cos(2\theta(t,E))   
\right] \, dt \right|   
\\ 
\leq C n^{-1/2}  
\end{multline*}  
for $\eta \in [n^3, (n+1)^3]$   
by the decay of $q$ and the estimate \eqref{aaa1}. Putting this into   
\eqref{longint}, we bound the last term in the expression   
for $\Theta(\xi,E)$ by $Cn^{-1}$.   
\end{proof}   
  
Next we need Lemma 8.4 from \cite{KLS}, concerning series   
of random variables. We will denote by $\rm{Exp}(X)$ the expectation   
of the random variable $X.$  
Assume $X_i(\omega)$ are independent random variables  
with zero mean ${ \rm Exp}(X_i)=0,$  and let   
\begin{equation}  
\label{rantyp}  
Z_n = X_n f_n(X_1, \dots, X_{n-1}),  
\end{equation}  
 where $f_n$ are some measurable   
functions. We have   
\begin{lemma}[\cite{KLS}]  
\label{ranser}  
Suppose that ${\rm Exp}(Z_n^2)\leq Cn^{-2\alpha}.$ Then for a.e.~$\omega,$ \\  
{\rm (1)} If $\alpha<\frac12$ and $\beta >\frac12 (1-2\alpha),$ then   
\[ \lim\limits_{n \rightarrow \infty} \left| \sum\limits_{j=1}^n  
 Z_j \right|   
n^{-\beta} =0. \]  
{\rm (2)} If $\alpha = \frac12$ and $\beta>\frac12,$ then   
\[ \lim\limits_{n \rightarrow \infty} \left| \sum\limits_{j=1}^n  
 Z_j \right|   
(\log n)^{-\beta} =0. \]  
{\rm (3)} If $\alpha>\frac12,$   
\[ \lim\limits_{n \rightarrow \infty} \left| \sum\limits_{j=1}^n  
 Z_j \right|   
=Y_\infty \]  
exists, and for any $\beta < \alpha -\frac12,$   
\[ \lim\limits_{n \rightarrow \infty} \left| \sum\limits_{j=n}^\infty  
 Z_j \right|   
n^{\beta} =0. \]  
\end{lemma}  
  
Consider   
\begin{multline}  
\label{conint1}  
 \log \frac{R((n+1)^3,E)}{R(n^3,E)} = I_n (E, \omega)  
\\ 
=   
\frac{1}{2n^{1/2}} \int\limits_{n^3}^{(n+1)^3}    
\frac{f(\xi^{1/3}-n)}{\xi^{2/3}}  
 \sin(2\xi +a_n(\omega)) \sin 2\theta(\xi,E)\,d\xi.   
\end{multline}  
We have   
\begin{lemma}  
\label{keyranle}  
The expression $I_n(E,\omega)$ may be represented as   
\begin{equation}  
\label{ranrep1}  
I_n(E, \omega) = \frac{9\pi}{8n} |\hat{f}(\frac{3E}{c})|^2 +   
Z_n^{(1)}+Z_n^{(2)} 
\end{equation}  
where $|Z_n^{(1)}(\omega)| \leq Cn^{-3/2}$, $\{Z_n^{(2)}(\omega)\}$   
is a sequence of random variables of type \eqref{rantyp}, 
and  $(Z_n^{(2)})^2$ has expectation $\leq Cn^{-1}$. 
\end{lemma}  
 
\begin{proof}  
The first observation is that we can replace 
$2\theta(\xi,E)$ in \eqref{conint1} by   
\[ 2\xi + 2F(\xi,E)+\theta_n+  
\int\limits_{n^3}^\xi \frac{q(c\eta^{2/3})}{c \eta^{2/3}}  
\cos (2\eta+ 2F(\eta,E) + 2\theta_n)\,d\eta. \]  
By Lemma~\ref{thetarep}, the value in \eqref{conint1} will change by at 
most   $Cn^{-3/2},$ and we can put this difference into $Z_n^{(1)}.$   
Also by an estimate similar to \eqref{devest},   
\begin{eqnarray}   
\nonumber \sin \left(  
2\xi + 2F(\xi,E)+ 2\theta_n+  
\int\limits_{n^3}^\xi \frac{q(c\eta^{2/3})}{c \eta^{2/3}}  
\cos (2\eta+ 2F(\eta,E) + 2\theta_n)\,d\eta  
 \right)   
\\   
= 
\label{terr} \cos (2\xi + 2F(\xi,E)+2\theta_n)  
 \int\limits_{n^3}^\xi \frac{q(c\eta^{2/3})}{c \eta^{2/3}}  
\cos (2\eta+ 2F(\eta,E) + 2\theta_n)\,d\eta   
\\ \nonumber    
+ 
\sin \left( 2\xi + 2F(\xi,E)+2\theta_n \right) + O(n^{-1}),   
\end{eqnarray}  
and so it suffices to consider the contributions from the   
first two terms on the right hand side,   
putting the $O(n^{-1})$ terms (multiplied by $n^{-1/2}$) into $Z_n^{(1)}.$  
We will write $Z_n^{(2)}$ as a sum of four parts, 
each satisfying the conclusion of the lemma. 
The first of these four is 
\[  
\frac{1}{2n^{1/2}} \int\limits_{n^3}^{(n+1)^3} \xi^{-2/3}   
f(\xi^{1/3}-n) \sin(2\xi +a_n(\omega)) \sin \left( 2\xi + 2F(\xi,E)+  
2\theta_n\right)\,d\xi. \]  
Clearly its mean is zero, and the expectation of its  
square is $O(n^{-2})$. It remains to discuss the   
contribution arising from the  first term on the right hand   
side in \eqref{terr}:  
\begin{multline} 
\label{lastt}  
\frac{1}{2n} \int\limits_{n^3}^{(n+1)^3} \xi^{-2/3}   
f(\xi^{1/3}-n) \sin(2\xi +a_n(\omega)) \cos (2\xi + 2F(\xi,E)+2\theta_n)  
\\  
\qquad\qquad \times  
\int\limits_{n^3}^\xi  \eta^{-2/3}   
f(\eta^{1/3}-n) \sin(2\eta +a_n(\omega))  
\cos (2\eta+ 2F(\eta,E) + 2\theta_n)\,d\eta d\xi   
\\  
= 
\frac{1}{4n} \left( \int\limits_{n^3}^{(n+1)^3} \xi^{-2/3}   
f(\xi^{1/3}-n) \sin(2\xi +a_n(\omega)) \cos (2\xi + 2F(\xi,E)+2\theta_n)  
\right)^2.  
\end{multline}  
Write the product of trigonometric functions as    
\begin{multline*}  
 \sin(2\xi +a_n(\omega)) \cos (2\xi + 2F(\xi,E)+2\theta_n)  \\  
= 
\frac{1}{2} \left(   
\sin(4\xi +2F(\xi,E)+2\theta_n+  
a_n(\omega))+ \sin ( 2F(\xi,E)+2\theta_n-a_n(\omega)) \right).  
\end{multline*}  
The first term, by integration by parts, gives a contribution of   
the order $n^{-3}$.  
Introducing notation   
\begin{align*}  
 I_n^{\sin} &= \int\limits_{n^3}^{(n+1)^3} \xi^{-2/3} f(\xi^{1/3}-n)   
\sin(2F(\xi,E) + 2\theta_n)\,d\xi, \\  
 I_n^{\cos} &= \int\limits_{n^3}^{(n+1)^3} \xi^{-2/3} f(\xi^{1/3}-n)   
\cos(2F(\xi,E) + 2\theta_n)\,d\xi  
\end{align*}  
we obtain that, modulo small corrections that can be put into $Z_n^{(1)},$   
the right-hand side of \eqref{lastt} is equal to  
\begin{equation}  
\label{termleft}   
 \frac{1}{8n}\left(  (I_n^{\sin})^2 (\cos(a_n(\omega)))^2 +  
 (I_n^{\cos})^2 (\sin(a_n(\omega)))^2 + I_n^{\sin}I_n^{\cos}\sin(2a_n(\omega))  
\right).   
\end{equation}  
The contribution of the last term in the brackets is one of the random   
variables of type \eqref{rantyp} constituting $Z_n^{(2)}.$ Subtracting   
from the first two terms their expectations, we get the other two constituents 
of $Z_n^{(2)}$.  
Finally, note that ${\rm Exp}((\cos(a_n(\omega)))^2)=  
{\rm Exp}((\sin(a_n(\omega)))^2)=\pi,$ and so after subtracting $Z_n^{(2)}$  
from \eqref{termleft} and taking expectations we deduce that the right-hand 
side of \eqref{lastt}, modulo the sum of the various terms $Z_n$, equals  
\[ \frac{\pi}{8n}\left(  (I_n^{\sin})^2  +  
 (I_n^{\cos})^2 \right) = \frac{\pi}{8n} \Big|  
\int\limits_{n^3}^{(n+1)^3} \xi^{-2/3} f(\xi^{1/3}-n) e^{2iF(\xi,E)}\,d\xi  
\Big|^2. \]  
Recalling that $F(\xi,E) = \frac{3}{2c}E(\xi^{1/3}-n),$ and making a change  
of variable $t=\xi^{1/3},$ we obtain  
\[   \frac{9\pi}{8n} \Big| \hat{f}({3E}/{c}) \Big|^2\ , \]  
proving the lemma.  
\end{proof}  
  
\begin{proof}[Proof of Theorem~\ref{rangrowth}]  
The proof follows directly from combining results of 
Lemma~\ref{keyranle},   
Lemma~\ref{ranser}, and the relation \eqref{intsame},   
choosing $\beta<1$ in conclusion (2) of Lemma~\ref{ranser}.
\end{proof}  
  
The final step in our construction is showing that for a.e.~$\omega,$   
a.e.~$E$ there exists a subordinate solution.   
\begin{proposition}  
\label{decsolp}  
Define the potential $q$ by \eqref{ranpot}. Consider any 
energy $E$ such that $\hat{f}(\frac{3E}{c}) \ne 0.$ Then for a.e.~$\omega$  
there exists a boundary condition $\beta(\omega)= \beta(\omega,E)$ such that   
\begin{equation}  
\label{decsol}  
\lim\limits_{\xi \rightarrow \infty} \frac{\log R(\xi,E,\beta(\omega))}  
{\log \xi} = -\Lambda (E),   
\end{equation}  
where $\Lambda(E)= (3\pi/8)\,|\hat f(3E/c)|^2$. 
\end{proposition}    
\begin{proof}  
By Lemma 8.7 of \cite{KLS} (which is formulated there in a   
discrete setting but   
extends trivially to our situation), it is sufficient to show that   
for a.e.~$\omega,$   
\[ \rho(\xi,E) = \log \frac{R(\xi,E,0)}{R(\xi,E,\pi/2)} \]  
has a limit $\rho_\infty(E)$ as $\xi \rightarrow \infty,$ and that   
\begin{equation}  
\label{rhoest}  
 \limsup\limits_{\xi \rightarrow \infty} \frac{\log |\rho(\xi,E)-\rho_\infty(E)|}  
{\log \xi} \leq -2\Lambda(E).   
\end{equation}  
Notice that by constancy of the Wronskian,  
\[ R(\xi,E,0)R(\xi,E,\frac{\pi}{2}) \sin(\theta(\xi,E,0)-  
\theta(\xi,E,\pi/2)) = \text{constant.} \]  
Therefore, by Theorem~\ref{rangrowth} for a.e.~$\omega$   
\begin{equation}  
\label{shut}  
\lim\limits_{\xi \rightarrow \infty} \frac{\log |\theta(\xi,E,0)-  
\theta(\xi,E,\pi/2)|}{\log \xi} = -2\Lambda(E)  
\end{equation}  
On the other hand, by Lemma~\ref{thetarep},   
\begin{multline*}  
 \rho((n+1)^3,E) = \sum\limits_{j=1}^{n+1} n^{-1/2}   
\int\limits_{n^3}^{(n+1)^3} \xi^{-2/3} f(\xi^{1/3}-n) \sin (2\xi+a_n(\omega))  
\\  
\times  
\left( \sin(2\xi+2F(\xi,E)+2\theta_{n}(\xi,E,0))-  
 \sin(2\xi+2F(\xi,E)+2\theta_{n}(\xi,E,\pi/2)) \right)  
\end{multline*}   
plus $ O(n^{-1})$.  
Expanding $\sin (2\xi+a_n(\omega))$ into sum of products, and using   
\eqref{shut}, we can apply  
Lemma~\ref{ranser} to prove convergence to $\rho_\infty$   
and estimate the rate of convergence obtaining \eqref{rhoest}.   
The estimate \eqref{devest} as before allows us to pass from estimates   
over the sequence $n^3$ to estimates over all $\xi$.  
\end{proof}  
  
Our main Theorem~\ref{count} follows from Proposition~\ref{decsolp},  
 Theorem~\ref{rangrowth} and Lemma~\ref{derivest}  
immediately.   
\begin{proof}[Proof of Theorem~\ref{count}]  
Going back to the $x$ variable representation, we find that   
for a.e.~$\omega,$ for a.e.~$E$ there exist both a decaying solution   
$u_d(x,E)$ and a growing solution $u_g(x,E)$, and the asymptotic behavior   
of their $L^2$ norms taken over $[0,L]$ is given by    
\begin{eqnarray}    
\label{decas}  
\log \int\limits_0^L |u_d (x,E)|^2\,dx = {\rm max}(0,\frac{1}{2}-3  
\Lambda(E))\log L (1+o(1)), \\  
\label{groas}  
\log \int\limits_0^L |u_g (x,E)|^2\,dx = (\frac{1}{2}+3  
\Lambda(E)) \log L (1+o(1)).   
\end{eqnarray}   
Hence, given a generic $\omega,$  there is a solution subordinate on the right   
for a.e.~$E.$   
\end{proof}   
  
We would like to conclude by making several remarks about some curious  
properties of the random potential \eqref{qx}.   
First of all, the structure of  potential \eqref{qx}  
indicates that to get singular or point spectrum   
``at minimal cost'' in terms of decay or smoothness of random  
potential, it is important to have correlations over large (increasing)  
distances. In particular, in our example correlation distance grows  
as a square root of the distance from the origin.    
This contrasts with the free Schr\"odinger operator case \cite{KU,KLS}, where long   
distance correlation would have been ineffectual. Another characteristic   
feature is the presence of increasingly fast oscillations at infinity.   
  
The half-line Stark operator with potential \eqref{qx} exhibits a number  
of interesting spectral and dynamical properties which are similar  
to the properties of potentials with the borderline decay rate in a model of \cite{KLS}  
in the free case. Introduce notation $X$ for the coordinate operator acting  
by $Xf(x) = |x|f(x).$  
Denote also   
$\langle \psi_1, \psi_2 \rangle$ the inner product of two vectors.  
Given an initial state $\psi,$  
$\psi (t)$ stands for its unitary evolution $e^{-iH_q t} \psi.$  
One reasonable   
measure of propagation properties is the averaged moments of coordinate  
operator   
\[ \langle \langle X^m \psi(t), \psi(t) \rangle \rangle_T =   
\frac{1}{T} \int\limits_0^T  \langle \psi(t), X^m \psi(t) \rangle \, dt. \]  
Scaling  
\[  \langle \langle X^m \psi(t), \psi(t) \rangle \rangle_T \approx CT^m \]  
corresponds to the constant velocity ballistic rate, while the behavior  
\[  \langle \langle X^m \psi(t), \psi(t) \rangle \rangle_T \approx CT^{2m} \]  
to the superballistic constant  
acceleration rate, as for the free Stark operator.    
We will denote by $P_c$ the spectral projector on the continuous part of the   
spectrum of the operator $H_q.$ Also, given a Borel measure $\mu,$ we   
will say that $\mu$ has 
local Hausdorff dimension $d(E)$ at the point $E$   
if for every $\epsilon >0$ there exists 
$\delta>0$ such that the restriction   
of $\mu$ to $[E-\delta, E+\delta]$ gives zero weight to any set of   
Hausdorff dimension $d(E)-\epsilon,$ and is supported on a set of  
Hausdorff dimension $d(E)+\delta.$ Let us denote by $H_{q,\beta}$  
the Schr\"odinger operator \eqref{stark} defined on the positive half-line 
with the boundary condition  
$\beta u(0)-u'(0)=0.$ Finally, let $\mu_\beta$ be  
the spectral measure associated with  
$H_{q,\beta}$ in a canonical way, namely 
\[ \langle (H_{q,\beta}-z)^{-1} \delta_0, \delta_0 \rangle =  
\int \frac{d\mu_\beta(t)}{t-z}, \] 
where $\delta_0$ is the delta-function distribution at $0$  
 (see, e.g. \cite{Simvan}, Sections I.1-I.6).  
The following proposition holds.  
\begin{proposition}  
\label{exodyn}  
Assume that $q_\omega(x)$ is given by \eqref{qx}. 
Then for a.e.~$\beta,$ $\omega$ the spectral   
measure of $H_{q,\beta}$ is pure point on the set where $\Lambda (E) >1/6,$   
and is singular continuous with local Hausdorff dimension   
\[ d(E) = 1- 6\Lambda(E)  \]  
on the set where $\Lambda(E)<1/6.$    
Moreover, for a.e.~$\beta,$ $\omega,$  
for any $\psi$ such that $P_c(\beta,\omega) \psi \ne 0,$   
we have   
\[  \langle \langle X^m \psi(t), \psi(t) \rangle \rangle_T \geq C_{\epsilon,   
m, \beta,\omega} T^{2m(1-\epsilon)} \]  
for any $\epsilon >0.$   
\end{proposition}  
\begin{proof}  
The proof of this proposition is analogous to the proof of Theorem 8.10 of   
\cite{KLS} and Theorem 5.1   
of \cite{KL}, where similar properties were established for a model  
\eqref{KLSR} of \cite{KLS} with $\alpha = 1/2$ in a discrete setting.  
Let us sketch the arguments for the sake of completeness.   
First, let  
\[ A = \left\{ E \left| \Lambda(E)> 1/6 \right. \right\}. \] 
Notice that by \eqref{lyapexp} and properties of $f,$ $A$ is a finite union  
of disjoint bounded intervals. By Proposition~\ref{decsolp}, for  
a.e.~$\omega,$ for a.e.~$E \in A,$ there exists a solution  
$u_{\beta(\omega)}(x,E,\omega) \in L^2.$ Also, the set $A$ belongs to the  
spectrum of $H_{q,\beta}$ for any $\omega$ and $\beta$ since $q$ is decaying. 
The general theory of rank one perturbations then implies that for  
a.e.~$\omega,$ for a.e.~$\beta$ the spectrum on $A$ is dense pure point  
(see, e.g. \cite{DRMS}, Theorem 5.1, or \cite{Siman,Simvan}). 
 
Now let $E$ be such that $\Lambda(E)<1/6,$ and fix $\epsilon>0.$  
Consider an interval $I_\delta = [E-\delta, E+\delta]$ such that  
the values of the Hausdorff dimension function $d(E)=1-6\Lambda(E)$ on  
$I_\delta$ belong to
the interval $(d(E)-\epsilon,d(E)+\epsilon)$. 
We can  choose $\delta>0$ since $\Lambda(E)$ is continuous.  
Denote   
\[ \|u\|^2_L = \int\limits_0^L |u(x)|^2\,dx, \]  
and recall that 
the upper $\alpha$ derivative of the measure $\mu$ is defined by   
\[ D^\alpha \mu(E) = \limsup_{\delta \rightarrow 0}  
\frac{\mu(E-\delta, E+\delta)}  
{(2\delta)^\alpha}. \]  
By the Jitomirskaya-Last  
extension of the subordinacy theory \cite{JL},  
the spectral measure $\mu_\beta$  
satisfies   
\begin{equation}  
\label{jl}  
 D^\alpha \mu_\beta(E) =  
\infty \,\,\,{\rm iff}\,\,\, \liminf_{L \rightarrow \infty}  
\frac{\|u_\beta \|_L^{2-\alpha}}{\|u_{\beta^{\perp}} \|_L^\alpha}  
=0,   
\end{equation}  
where $u_\beta$ is a nonzero solution satisfying the boundary  
condition at zero,   
and $u_{\beta^\perp}$ is the solution satifying an orthogonal boundary  
condition $\beta^{-1}u(0)+u'(0)=0$ (in fact, any linearly independent  
solution will do in place of $u_{\beta^\perp}$).    
From the estimates \eqref{decas}, \eqref{groas}  
it follows that for a.e.~$\omega$ and a.e.~$E \in I_\delta,$ there  
exists a boundary condition $\beta(\omega)$ such that  
\begin{equation} 
\label{LHD} 
 {\rm inf} \left\{ \alpha \left| D^\alpha \mu_{\beta(\omega)}(E)=\infty 
\right. \right\}=d(E).  
\end{equation} 
Recall two basic facts from the general rank one perturbation theory: 
first, 
\[ \int_\reals\mu_\beta (S) d\beta = m(S) \] 
for any Borel measurable set $S,$  
where $m$ is the Lebesgue measure. Second,  
we can have $D^1 \mu_{\beta_1}(E)=D^1 \mu_{\beta_2}(E)=\infty$ for some  
$\beta_1 \ne \beta_2$ only on a fixed (for a given $\omega$) exceptional  
set of energies of Lebesgue measure zero.  
The first fact, 
in the context of Schr\"odinger operators,  
is due to Simon and Wolff \cite{SW}. The second is a simple consequence 
of an explicit formula relating Cauchy transforms of $\mu_{\beta_1}$ and 
$\mu_{\beta_2}.$ 
See, for example, equation (I.13) of \cite{Simvan} for more details.  
Using these properties, we see that \eqref{LHD} 
implies that for a.e.~$\omega,$ 
and a.e.~$\beta$ we have  
\[  {\rm inf} \left\{ \alpha \left| D^\alpha \mu_{\beta}(E)=\infty 
\right. \right\}=d(E) \] 
for every $E \in I_\delta$ in the support of $\mu_\beta.$  
Similarly, we also obtain that for a.e.~$\omega,$ $\beta$ the solution  
$u_\beta(x,E)$ satisfies the bound \eqref{decas} for the energies in the  
support of $\mu_\beta.$   
By the choice of $\delta$ and well-known relation between $D^\alpha$  
derivatives and dimensional properties of Radon measures (see \cite{Ro})  
it follows that for a.e.~$\omega,$ $\beta$ the restriction of the  
spectral measure $\mu_\beta$ to $I_\delta$ is supported on a set of  
Hausdorff dimension $\leq d(E)+\epsilon$ and gives zero weight to any set  
of Hausdorff dimension $\leq d(E)-\epsilon.$ 
   
To show the dynamical bound, notice that Theorem 1.2 of \cite{KL} says 
 that if the   
spectral measure $\mu$ is $\alpha-$continuous on a set $S,$ and for 
 every $E \in S,$   
the generalized eigenfunction $u(x,E)$ satisfies   
\[ \limsup_{L\rightarrow \infty} L^{-\gamma} \|u\|^2_L <\infty, \]  
then for any vector $\psi$ with nonzero projection on $S,$ 
that is, with $P_S (\psi) \ne 0,$   
\[ \langle \langle |X|^m \psi(t), \psi(t) \rangle  \rangle_T 
 \geq CT^{\frac{m \alpha}{\gamma}}. \]  
The proof is now completed as in Theorem 5.1 in \cite{KL},   
noticing that in our context, $\alpha \sim 1-6\Lambda(E),$ while by 
 \eqref{decas}  
$\gamma \sim \frac{1}{2} - 3 \Lambda(E).$    
\end{proof}  
 
\noindent \it Remarks \rm 1. The methods of \cite{KL} also give 
a stronger dynamical estimate that for a.e.~$\beta,$ $\omega$ and for 
every $\epsilon>0$ and $\rho>0,$ there exists a constant $C_{\beta, 
\omega,\rho,\epsilon}$ such that if $R_T = C_{\beta, 
\omega,\rho,\epsilon}T^{2-\epsilon},$ then  
\begin{equation}\label{DE} 
\langle \|\psi(t)\|_{R_T}^2 \rangle_T \leq \|\psi - P_c(\beta,\omega)\psi\|^2 
+\rho \|\psi\|^2. 
\end{equation} 
From \eqref{DE} it follows that the whole part of the 
wavepacket lying in the continuous spectral subspace travels at  a rate  
$\gtrsim T^{2(1-\epsilon)},$ 
even though we can choose $f$ in \eqref{qx} so that 
the spectral dimension $d(E)$ is arbitrarily close to zero in some parts  
of the support of the spectral measure. \\ 
2. The Proposition~\ref{exodyn} is not true for a fixed boundary 
condition. In fact, for a.e.~$\omega,$  
for a dense $G_\delta$ set of boundary conditions $\beta \in \reals,$ 
the spectrum of $H_{q,\beta}$ is going to be purely singular continuous  
\cite{DRMS}. \\ 
 
A final remark, which is more of academic interest, is that our method   
can be modified in a straightforward way to treat perturbations of   
background operators   
\[ -\frac{d^2}{dx^2} - x^{\lambda}, \]  
where   
$0<\lambda<2$ (the operator becomes non-selfadjoint for larger $\lambda$).   
In particular, one obtains that the absolutely continuous spectrum of such   
operators is preserved for perturbations $q$ satisfying   
$|q(x)| \leq C(1+|x|)^{-\frac{1}{4}(2-\lambda)-\epsilon}$ or   
$ q \in C^{\frac{1}{2\lambda}(2-\lambda)+\epsilon}(\reals),$ where  
$\epsilon$ is an arbitrary small positive number. These results are  
optimal in the power scale and in H\"older spaces, respectively.

\end{document}